\documentstyle[twoside]{article}
\title{Asymptotics of the Wright function ${}_1\Psi_1(z)$ on the Stokes lines}
\author{\sc R. B.\ Paris \\
{\em School of Engineering, Computing and Applied Mathematics}, \\
{\em University of Abertay Dundee, Dundee DD1 1HG, UK}}
\begin{document}
\def\f#1#2{\mbox{${\textstyle \frac{#1}{#2}}$}}
\def\dfrac#1#2{\displaystyle{\frac{#1}{#2}}}
\def\boldal{\mbox{\boldmath $\alpha$}}
{\newcommand{\Sgoth}{S\;\!\!\!\!\!/}
\newcommand{\bee}{\begin{equation}}
\newcommand{\ee}{\end{equation}}
\newcommand{\lam}{\lambda}
\newcommand{\ka}{\kappa}
\newcommand{\al}{\alpha}
\newcommand{\th}{\theta}
\newcommand{\fr}{\frac{1}{2}}
\newcommand{\fs}{\f{1}{2}}
\newcommand{\g}{\Gamma}
\newcommand{\br}{\biggr}
\newcommand{\bl}{\biggl}
\newcommand{\ra}{\rightarrow}
\newcommand{\mbint}{\frac{1}{2\pi i}\int_{c-\infty i}^{c+\infty i}}
\newcommand{\mbcint}{\frac{1}{2\pi i}\int_C}
\newcommand{\mboint}{\frac{1}{2\pi i}\int_{-\infty i}^{\infty i}}
\newcommand{\gtwid}{\raisebox{-.8ex}{\mbox{$\stackrel{\textstyle >}{\sim}$}}}
\newcommand{\ltwid}{\raisebox{-.8ex}{\mbox{$\stackrel{\textstyle <}{\sim}$}}}
\renewcommand{\topfraction}{0.9}
\renewcommand{\bottomfraction}{0.9}
\renewcommand{\textfraction}{0.05}
\newcommand{\mcol}{\multicolumn}
\date{}
\maketitle
\pagestyle{myheadings}
\markboth{\hfill \sc R. B.\ Paris  \hfill}
{\hfill \sc  Wright function asymptotics\hfill}
\begin{abstract}
We investigate a particular aspect of the asymptotic expansion of the 
Wright function ${}_1\Psi_1(z)$ for large $|z|$. The form of the exponentially small expansion associated with this function on certain rays in the $z$-plane (known as Stokes lines) is discussed. The main thrust of the paper is concerned with the expansion in the particular case when the Stokes line coincides with the negative real axis $\arg\,z=\pi$.
Some numerical examples which confirm the accuracy of the expansion are given.
\vspace{0.4cm}

\noindent {\bf Mathematics Subject Classification:} 33C20, 33C70, 34E05, 41A60 
\vspace{0.3cm}

\noindent {\bf Keywords:} Asymptotics, exponentially small expansions, Wright function, Stokes lines
\end{abstract}

\vspace{0.3cm}

\noindent $\,$\hrulefill $\,$

\vspace{0.2cm}

\begin{center}
{\bf 1. \  Introduction}
\end{center}
\setcounter{section}{1}
\setcounter{equation}{0}
\renewcommand{\theequation}{\arabic{section}.\arabic{equation}}
We consider the subclass of Wright functions defined by
\bee\label{e11}
{}_1\Psi_1(z)\equiv{}_1\Psi_1\left(\begin{array}{c}a, \al\\ b, \beta \end{array};z\right)=\sum_{r=0}^\infty \frac{\g(\al r+a)}{\g(\beta r+b)}\,\frac{z^r}{r!}, 
\ee
where the parameters $\alpha$ and $\beta$ are real and positive and $a$ and $b$ are
arbitrary complex numbers. We also assume throughout that the $\alpha$ and $a$ are subject to the restriction
\bee\label{e11a}
\alpha r+a\neq 0, -1, -2, \ldots \qquad (r=0, 1, 2, \ldots)
\ee
so that no gamma function in the numerator in (\ref{e11}) is singular.
This function has recently found application in probability theory concerned with refined local approximations for members of some Poisson-Tweedie exponential dispersion models \cite{PV}.
In the special case $\alpha=\beta=1$, the function ${}_1\Psi_1(z)$ reduces to 
the confluent hypergeometric function $\g(a)/\g(b)\,{}_1F_1(a;b;z)$;
see, for example, \cite [ p.~40]{S}. When $\alpha=1$, the function ${}_1\Psi_1(z)$ becomes the generalised Mittag-Leffler function considered in \cite{Pr}.

We introduce the parameters given by
\bee\label{e12}
\kappa=1+\beta-\alpha, \qquad 
h=\alpha^{\alpha}\beta^{-\beta},\qquad\vartheta=a-b.
\ee
If it is supposed that $\alpha$ and $\beta$ are such that $\kappa>0$ then ${}_1\Psi_1(z)$ 
is uniformly and absolutely convergent for all finite $z$. If $\kappa=0$, the sum in (\ref{e11})
has a finite radius of convergence equal to $h^{-1}$, whereas for $\kappa<0$ the sum is divergent 
for all nonzero values of $z$. The parameter $\kappa$ is found to play a critical role 
in the asymptotic theory of ${}_1\Psi_1(z)$ by determining the sectors in the $z$ plane 
in which its behaviour is either exponentially large or algebraic in character as $|z|\ra\infty$. 

The asymptotic expansion of ${}_1\Psi_1(z)$ for $|z|\ra\infty$ and finite 
values of the parameters can be obtained from 
\cite{Br, W1, W2}; see also \cite[ \S 2.3]{PK}. In Section 2 we present a summary of the standard expansion of ${}_1\Psi_1(z)$ when the parameter $\kappa$ satisfies $0<\kappa<2$. In this case, the function has a composite expansion consisting of a single exponential expansion and an algebraic expansion. 
When $0<\kappa\leq1$, which is the case that principally concerns us here, a Stokes phenomenon occurs on the rays $\arg\,z=\pm\pi\kappa$, where the exponential expansion is maximally subdominant with respect to the algebraic expansion. A more precise understanding of the asymptotic behaviour of the Wright function is then achieved by taking  this phenomenon into account \cite{PLMJ}.

Our objective in this paper is to examine the form of the exponentially small expansion associated with ${}_1\Psi_1(z)$ in the special case when the parameter $\kappa=1$, where the Stokes line coincides with the negative real axis $\arg\,z=\pi$. The approach we adopt to determine this expansion is by means of a Mellin-Barnes integral representation for ${}_1\Psi_1(z)$ combined with the theory presented in \cite[Chapter 6]{PK}.
Numerical verifications confirming the validity of this expansion are presented. 

\vspace{0.6cm}

\begin{center}
{\bf 2. \ Standard asymptotic theory of ${}_1\Psi_1(z)$ for $|z|\ra\infty$}
\end{center}
\setcounter{section}{2}
\setcounter{equation}{0}
\renewcommand{\theequation}{\arabic{section}.\arabic{equation}}
In this section we present the standard asymptotic expansion of ${}_1\Psi_1(z)$ as 
$|z|\ra\infty$ for $0<\kappa<2$ and finite values of the parameters given in \cite{Br} and 
\cite{W2}; see also \cite[ \S2.3]{PK}. Throughout we let $\epsilon$ denote an arbitrarily small positive quantity. We first define the associated exponential and algebraic expansions $E(z)$ and $H(z)$ given by the formal asymptotic sums
\bee\label{e22c}
E(z):=Z^\vartheta e^Z\sum_{j=0}^\infty A_jZ^{-j}, \qquad Z=\kappa (hz)^{1/\kappa},
\ee
and
\bee\label{e25}
H(z):=\frac{1}{\alpha}\sum_{k=0}^\infty \frac{(-)^k}{k!}\,\frac{\Gamma((k+a)/\alpha)}{\Gamma(b-\beta(k+a)/\alpha)} z^{-(k+a)/\alpha}.
\ee
The $A_j$ appear as coefficients in the inverse factorial expansion of the quotient of gamma functions given by \cite[\S 3]{Br}, \cite[p.~39]{PK}
\newtheorem{lemma}{Lemma}
\begin{lemma}
Let $M$ denote a positive integer and suppose that $\kappa>0$. Then there exist coefficients $A_j$ $(0\leq j\leq M-1)$ such that
\bee\label{e22b}
\g(s)\,\frac{\g(1-b+\beta s)}{\g(1-a+\alpha s)}
=\kappa(h\kappa^\kappa)^{-s}\left\{\sum_{j=0}^{M-1}(-)^jA_j\g(\kappa s+\vartheta-j)+\sigma_M(s) \g(\kappa s+\vartheta-M)\right\},
\ee
where the parameters $\kappa$, $h$ and $\vartheta$ have the values given in (\ref{e12}). The remainder function
$\sigma_M(s)$ is analytic in $s$ except at the poles of the corresponding gamma function ratio and is such that $\sigma_M(s)=O(1)$ as $|s|\ra\infty$ uniformly in $|\arg\,s|\leq\pi-\epsilon$.
\end{lemma}
The coefficients $A_j$ are independent of $s$ and depend only on the parameters $\alpha$, 
$\beta$, $a$ and $b$. The leading coefficients $A_0$ and $A_1$ are specified by \cite{P10}
\bee\label{e23}
A_0=\kappa^{-\frac{1}{2}-\vartheta}
\alpha^{a-\frac{1}{2}}\beta^{\frac{1}{2}-b},\qquad
A_1=\fs\kappa A_0\left(\frac{a^2-a+\f{1}{6}}{\alpha}-\frac{b^2-b+\f{1}{6}}{\beta}+
\frac{1-\kappa-6\vartheta(1-\vartheta)}{6\kappa}\right).
\ee
An algorithm for the numerical evaluation of the $A_j$ for a general quotient of gamma functions has been described in \cite{P10}; see the appendix for the specific case in (\ref{e22b}).

The expansion theorem for ${}_1\Psi_1(z)$ is then given by the following theorem.
\newtheorem{theorem}{Theorem}
\begin{theorem}
If $0<\kappa<2$, then
\bee\label{e21}{}_1\Psi_1(z)\sim\left\{\begin{array}{lll}E(z)+H(ze^{\mp\pi i}) & 
\mbox{in} & |\arg\,z|\leq \f{1}{2}\pi\kappa \\
\\H(ze^{\mp\pi i}) & 
\mbox{in} & |\arg (-z)|\leq\f{1}{2}\pi(2-\kappa)-\epsilon\end{array} \right.
\ee
as $|z|\ra\infty$. The upper or lower sign 
in $H(ze^{\mp\pi i})$ is chosen according as $z$ lies in the upper or lower 
half-plane, respectively.
\end{theorem}
This result gives the dominant expansion of ${}_1\Psi_1(z)$. It is seen that the $z$-plane is 
divided into two sectors, with a common vertex at $z=0$, by the rays (the anti-Stokes lines) 
$\arg\,z=\pm\fs\pi\kappa$. In the sector $|\arg\,z|<\fs\pi\kappa$, 
the asymptotic character of $_1\Psi_1(z)$ is exponentially large whereas in the complementary sector 
$|\arg (-z)|<\fs\pi(2-\kappa)$, the dominant expansion is algebraic in character. 

In the case  $0<\kappa<1$, 
the exponential expansion $E(z)$ is still present beyond the sector $|\arg\,z|<\fs\pi\kappa$ where it becomes subdominant in the sectors $\fs\pi\kappa<|\arg\,z|<\pi\kappa$. The rays $\arg\,z=\pm\pi\kappa$, where $E(z)$ is {\it maximally} subdominant with respect to $H(ze^{\mp\pi i})$, are called Stokes lines.\footnote{The positive real axis $\arg\,z=0$ is also a Stokes line where the algebraic expansion is maximally subdominant.} As these rays are crossed (in the sense of increasing $|\arg\,z|$) the exponential expansion  switches off according to the now familiar error-function smoothing law \cite{By}. In view of this interpretation of the Stokes phenomenon a more precise version of Theorem 1 in the case $0<\kappa<1$ was shown in \cite{PLMJ} to be given by the following.
\begin{theorem}
If $0<\kappa<1$, then 
\bee\label{e22}{}_1\Psi_1(z)\sim\left\{\begin{array}{lll}
E(z)+H(ze^{\mp\pi i}) & 
\mbox{in} & |\arg\,z|\leq \pi\kappa-\epsilon \\
\\ \fs {\cal E}(z)+H^o(ze^{\mp\pi i}) & \mbox{on} & \arg\,z=\pm\pi\kappa\\ 
\\H(ze^{\mp\pi i}) & 
\mbox{in} & |\arg (-z)|\leq\pi(1-\kappa)-\epsilon\end{array} \right.
\ee
as $|z|\ra\infty$. The upper or lower sign 
in $H(ze^{\mp\pi i})$ is chosen according as $z$ lies in the upper or lower 
half-plane, respectively. 

In the middle expression on the Stokes lines $\arg\,z=\pm\pi\kappa$, the superscript o denotes that the algebraic expansions $H(ze^{\mp\pi i})$ are optimally truncated at, or near, the term of least magnitude. 
The expansion ${\cal E}(z)$ denotes the expansion $E(z)$ augmented by the presence of an additional series and is given by
\[{\cal E}(z):=E(z)-\frac{2i}{\sqrt{2\pi X}}(Xe^{-\pi i})^\vartheta e^{-X}\sum_{j=0}^\infty (-)^jB_j X^{-j},\quad X=|Z|,\]
where $Z$ is defined in (\ref{e22c}) and the coefficients $B_j$ are defined in (\ref{e36}).
\end{theorem}
Although the expansion in (\ref{e21}) is a valid asymptotic description of ${}_1\Psi_1(z)$,
more accurate evaluation will result from taking the Stokes phenomenon into account as the Stokes rays are crossed. 
Theorem 1 deals with the expansion of ${}_1\Psi_1(z)$ only when $0<\kappa<2$.
When $\kappa\geq 2$, there are additional exponential series of the type $E(z)$ with suitably rotated arguments;  see, for example, \cite[p.~58]{PK} for details.

This concludes our summary of the expansion of ${}_1\Psi_1(z)$ when $0<\kappa<2$. We now confine our attention in the remainder of this paper to the case $\alpha=\beta$ (so that $\kappa=1$), which is not covered by Theorem 2, and consider the expansion of ${}_1\Psi_1(z)$ on the Stokes line  $\arg\,z=\pi$.

\vspace{0.6cm}

\begin{center}
{\bf 3. The expansion of ${}_1\Psi_1(z)$ on $\arg\,z=\pi$ when $\kappa=1$}
\end{center}
\setcounter{section}{3}
\setcounter{equation}{0}
\renewcommand{\theequation}{\arabic{section}.\arabic{equation}}
In this section we consider the particular Wright function with $\alpha=\beta$, viz.
\[{}_1\Psi_1\left(\begin{array}{c}a, \alpha\\b, \alpha\end{array};z\right)\equiv{}_1\Psi_1(z)\]
for which the parameters in (\ref{e12}) are $\kappa=h=1$. The Stokes line for this function is the negative real axis $\arg\,z=\pi$. Our aim is to derive the asymptotic expansion of ${}_1\Psi_1(-x)$ as $x\ra+\infty$, taking into account the exponentially small contribution. The detailed analysis of the more general Wright function ${}_1\Psi_q(z)$, with $q$ denominatorial gamma functions, has been given in \cite{PLMJ} when $0<\kappa<1$; here we present a summary for the case $q=1$ so that the paper is self-contained.

Our starting point is the Mellin-Barnes integral representation \cite[\S 2.4]{PK} 
\bee\label{e31}
{}_1\Psi_1(-x)=\frac{1}{2\pi i}\int_{-\infty i}^{\infty i} \g(s)\,\frac{\g(a-\alpha s)}{\g(b-\alpha s)}\,x^{-s}\,ds\qquad (|\arg\,x|<\fs\pi),
\ee
where the integration path is indented 
to separate the two sequences of poles of the integrand at $s=k$ and $s=(k+a)/\alpha$ ($k=0, 1, 2, \ldots $).
Displacement of the integration path to the right in the usual manner over the first $m$ poles of $\g(a-\alpha s)$  then produces
\bee\label{e32}
{}_1\Psi_1(-x)=\frac{1}{\alpha}\sum_{k=0}^{m-1} \frac{(-)^k}{k!}\,\frac{\g((k+a)/\alpha)}{\g(b-a-k)}x^{-(k+a)/\alpha}+P_m(z),
\ee
where, upon use of the reflection formula for the gamma function,
\bee\label{e33}
P_m(z)=\frac{1}{2\pi i} \int_{L_m} \frac{\g(s)\g(1-b+\alpha s)}{\g(1-a+\alpha s)}\,\frac{\sin \pi(\alpha s-b)}{\sin \pi(\alpha s-a)}\,x^{-s}ds.
\ee
The path $L_m$ denotes a path (possibly indented) parallel to the imaginary $s$-axis with $\Re(s)=(a+m-c)/\alpha$, $0<c<1$. We recognise that the series in (\ref{e32}) is the first $m$ terms of the algebraic expansion $H(x)$ in (\ref{e25}).

In order to detect the appearance of the exponential series as $|z|\ra\infty$ we need to 
choose $m$ in the algebraic expansion (\ref{e32}) to correspond to the optimal truncation index $m_o$; that is truncation at, or near, the least term in magnitude. Making use of the well-known result
$\g(k+a)/\g(k+b)\sim k^{a-b}$ for $k\ra+\infty$, we easily find that
\bee\label{e31b}
m_o\simeq \alpha x.
\ee
We now set $m=m_o$ in (\ref{e32}) and so consider the algebraic expansion $H(x)$ to be {\it optimally truncated}.

When $m$ is chosen as above, it follows from (\ref{e31b})  that $m_o\ra\infty$ as $x\ra\infty$, so that on the integration path $L_m$ in (\ref{e33}) we have $|s|$ everywhere large. Consequently, we may employ the inverse factorial expansion (\ref{e22b}) for the ratio of gamma functions to find from (\ref{e33}), with the variable $s$ replaced by $(s+a+m)/\alpha$,
\[P_m(x)=x^{-\mu(a+m)} \sum_{j=0}^{M-1}(-)^jA_j\bl\{e^{-\pi i\vartheta}I(xe^{\pi i\alpha})-e^{\pi i\vartheta} I(xe^{-\pi i\alpha})\br\}+R_{M,m}(x),\]
where
\[I(z):=\frac{\mu}{4\pi}
\int_{-c-\infty i}^{-c+\infty i} \g(\mu s+\nu-j)\,\frac{z^{-\mu s}}{\sin \pi s}ds\qquad(0<c<1;\ |\arg\,z|<\fs\pi+\pi\alpha)\]
and we have defined
\bee\label{e31c}
\nu:=\mu(a+m)+\vartheta,\qquad \mu:=1/\alpha.
\ee
The remainder term $R_{M,m}(x)$ is an integral involving the quantity $\sigma_M(s)$ in (\ref{e22b}). This is not discussed here since
it was established in \cite{PLMJ} that $R_{M,m}(x)=O(x^{\vartheta-M}e^{-x}$) as $x\ra+\infty$.

The integral $I(z)$ is now expressed in terms of the generalised terminant function introduced in \cite{PLMJ} by
\[T_\nu(\mu; z)=(ze^{-\pi i/\mu})^{-\nu} \exp\,[ze^{-\pi i/\mu}]\,\frac{\mu}{4\pi}\int_{-c-\infty i}^{-c+\infty i}\g(\mu s+\nu)\,\frac{z^{-\mu s}}{\sin \pi s}\,ds\]
when $|\arg\,z|<\fs\pi+\pi/\mu$ and $0<c<1$. In the case $\mu=1$ ($\alpha=1$), this reduces to the standard terminant function $T_\nu(z)\equiv T_\nu(1;z)$ expressed as a multiple of the incomplete gamma function $\g(1-\nu,z)$  
\[T_\nu(z)=\frac{e^{\pi i\nu}\g(\nu)}{2\pi i}\,\g(1-\nu,z)\] 
introduced in \cite{Olver}; see also \cite[p.~260]{PK}.
The rays $\arg\,z=\pm\pi/\mu$ are Stokes lines for $T_\nu(\mu;z)$.
The connection formula\footnote{When $\mu=1$ ($\alpha=1$) the connection formula (\ref{ecf}) reduces to that given in \cite[p.~260]{PK}.} satisfied by this function with the arguments $xe^{\pm\pi i\alpha}$ ($x>0$), is given by \cite[Appendix A]{PLMJ}
\bee\label{ecf}
\exp\,[x-xe^{-2\pi i/\mu}] \,T_\nu(\mu;xe^{-\pi i/\mu})=e^{2\pi i\nu/\mu}(T_\nu(\mu; xe^{\pi i/\mu})-1).
\ee
Then, after some straightforward algebra making use of (\ref{ecf}), we obtain
\bee\label{e35}
P_m(x)=x^\vartheta e^{-x}\sum_{j=0}^{M-1}(-)^jA_jx^{-j}\bl\{e^{-\pi i\vartheta}T_{\nu-j}(\mu; xe^{\pi i\alpha})-e^{\pi i\vartheta} (T_{\nu-j}(\mu; xe^{\pi i\alpha})-1)+O(x^{-M})\br\}.
\ee

With $m$ chosen according to (\ref{e31b}) it follows that the parameter $\nu\sim x$ as $x\ra+\infty$.
The expansion of the terminant function $T_{\nu-j}(\mu; xe^{\pi i/\mu})$ on its Stokes line when $\nu\sim x\ra+\infty$ is given by \cite[Appendix A]{PLMJ}
\bee\label{e301}
T_{\nu-j}(\mu; xe^{\pi i/\mu}) =  \frac{1}{2}-\frac{i}{\sqrt{2\pi x}}\left\{\sum_{k=0}^{N-1} g_{_{2k}}(\mu;j)\,(\fs)_{_k} (\fs x)^{-k}+O(x^{-N})\right\}
\ee
for $j=0, 1, 2, \ldots\,$, where $N$ is a positive integer and $(a)_k=\g(k+a)/\g(a)$ is the Pochhammer symbol. 
The coefficients $g_{2k}(\mu;j)$ appear in the expansion
\[\frac{\mu\tau^{\delta_j-1}}{(1-\tau^\mu)}\,\frac{d\tau}{dw}=-\frac{1}{w}+\sum_{k=0}^\infty g_k(\mu;j)w^k,\qquad \fs w^2=\tau-\log\,\tau-1,\quad 
\delta_j:=\nu-x-j,\]
where $\delta_j$ is bounded
and the first few even-order coefficients  $g_{2k}(\mu;j)\equiv 6^{-2k}{\hat g}_{2k}(\mu;j)$ are
\begin{eqnarray*}
{\hat g}_0(\mu;j)\!\!&=&\!\!\f{1}{6}-\delta_j+\fs\mu,\\
{\hat g}_2(\mu;j)\!\!&=&\!\!\f{1}{30}(2+45\mu+45\mu^2-90\delta_j(1+3\mu+\mu^2)+270\delta_j^2(1+\mu)
-180\delta_j^3), \\
{\hat g}_4(\mu;j)\!\!&=&\!\!\f{1}{140}(-65+105\mu+630\mu^2-210\mu^4-42\delta_j(5+90\mu+100\mu^2-6\mu^4)+1260\delta_j^2
(3\\
&&+10\mu+5\mu^2)-840\delta_j^3(10+15\mu+3\mu^2)+
1260\delta_j^4(5+3\mu)-1512\delta_j^5),\\
{\hat g}_6(\mu;j)\!\!&=&\!\!\f{1}{700}(7(-16-417\mu+225\mu^2-1008\mu^4+180\mu^6)-
6\delta_j(-973+1575\mu+9555\mu^2\\
&&-4410\mu^4+180\mu^6)+1890\delta_j^2(5+91\mu+112\mu^2-14\mu^4)-
1260\delta_j^3(91+336\mu\\
&&+210\mu^2-6\mu^4)+26460\delta_j^4(8+15\mu+5\mu^2)-22680\delta_j^5
(7+7\mu+\mu^2)\\
&&+7560\delta_j^6(7+3\mu)-6480\delta_j^7),\ldots .
\end{eqnarray*}

Substitution of (\ref{e301}) into (\ref{e35})  (where we put $M=N$
for convenience) then yields for $x\ra+\infty$
\bee\label{e302}
P_m(x)= x^\vartheta e^{-x}\bl\{\cos \pi\vartheta \sum_{j=0}^{M-1} (-)^jA_jx^{-j}-\frac{2\sin \pi\vartheta}{\sqrt{2\pi x}}\sum_{j=0}^{M-1} (-)^jB_jx^{-j}+O(x^{-M})\br\},
\ee
where the coefficients $B_j$ are given by
\bee\label{e36}
B_j=\sum_{k=0}^j(-2)^k (\fs)_{_k} A_{j-k}\,g_{2k}(\mu;j-k).
\ee
From (\ref{e32}) and (\ref{e35}) we then obtain the expansion given in the following theorem.
\begin{theorem}
Provided $\vartheta=a-b$ is non-integer, we have the asymptotic expansion
when $\kappa=1$
\bee\label{e38a}
{}_1\Psi_1(-x)-H^o(x)\sim x^\vartheta e^{-x}\bl\{\cos \pi\vartheta \sum_{j=0}^{\infty} (-)^jA_jx^{-j}-\frac{2\sin \pi\vartheta}{\sqrt{2\pi x}}\sum_{j=0}^{\infty} (-)^jB_jx^{-j}\br\}
\ee
as $x\ra+\infty$, where the algebraic expansion $H(x)$ is defined in (\ref{e25}) and the superscript o signifies that this series is optimally truncated. 
The coefficients $A_j$ are specified by the inverse factorial expansion (\ref{e22b})
and the coefficients $B_j$ are defined in (\ref{e36}).
\end{theorem}
\vspace{0.6cm}

\begin{center}
{\bf 4. The expansion of ${}_1\Psi_1(-x)$ when $\vartheta$ is an integer}
\end{center}
\setcounter{section}{4}
\setcounter{equation}{0}
\renewcommand{\theequation}{\arabic{section}.\arabic{equation}}
When $\vartheta=-n$, $n=1, 2, \ldots$,
the integrand in (\ref{e31}) has a finite set of poles at $s=(k+a)/\alpha$, $0\leq k\leq n-1$. Displacement of the integration path to the right over these poles produces
\bee\label{e37}
{}_1\Psi_1(-x)=
\frac{1}{\alpha}\sum_{k=0}^{n-1} \frac{(-)^k}{k!}\,\frac{\g((k+a)/\alpha)}{\g(b-a-k)}x^{-(k+a)/\alpha}
+\frac{(-)^n}{2\pi i}\int_{L}\frac{\g(s) \g(1-b+\alpha s)}{\g(1-a+\alpha s)}\,x^{-s}ds,
\ee
where $L$ denotes a path parallel to the imaginary $s$-axis with Re\,$(s)>{\mbox Re}\,(a)+n-1$. The above algebraic expansion contains $n$ terms and so cannot be optimally truncated as $x\ra+\infty$; the parameter $\nu$ in (\ref{e31c}) (with $m=n$) is therefore finite in this limit. As a consequence, the analysis in Section 3 based on the asymptotic structure of the terminant function of large order and argument is inapplicable in this case. 

We can now displace $L$ as far to the right as we please so that $|s|$ is everywhere large on the contour and the inverse factorial expansion (\ref{e22b}) may again be employed to produce the exponentially small expansion on the right-hand side of (\ref{e35}). When $\vartheta=-n$, the two terminant functions in (\ref{e35}) cancel to yield
\begin{eqnarray}
\frac{(-)^n}{2\pi i}\int_L\frac{\g(s) \g(1-b+\alpha s)}{\g(1-a+\alpha s)}\,x^{-s}ds\!\!
&=&\!\! x^\vartheta e^{-x}\cos \pi\vartheta \bl\{\sum_{j=0}^{M-1}(-)^jA_j x^{-j}+O(x^{-M})\br\},\label{e37a}
\end{eqnarray}
where we have  
put $(-)^n\equiv\cos \pi\vartheta$. 
An alternative procedure is that described in \cite{P10} which makes use of the Cahen-Mellin integral for $e^{-x}$.

From (\ref{e37}) and (\ref{e37a}) we then obtain the following theorem.
\begin{theorem}
When $\vartheta=a-b=-n$, $n=1, 2, \ldots\,$, we have the asymptotic expansion when $\kappa=1$
\bee\label{e38}
{}_1\Psi_1(-x)-\frac{1}{\alpha}\sum_{k=0}^{n-1} \frac{(-)^k}{k!}\,\frac{\g((k+a)/\alpha)}{\g(b-a-k)}x^{-(k+a)/\alpha}
\sim x^\vartheta e^{-x}\cos \pi\vartheta \sum_{j=0}^{\infty} (-)^jA_jx^{-j}
\ee
as $x\ra+\infty$, where the coefficients $A_j$ are specified by the inverse factorial expansion (\ref{e22b}).
\end{theorem}

In general the series on the right-hand side of (\ref{e38}) is an asymptotic series consisting of an infinite number of terms. It is of interest to note, however,  that when 
\bee\label{e4c}
\alpha=\frac{1}{q},\qquad a=\frac{p}{q},
\ee
where $p$ and $q$ denote relatively prime positive integers, it is found that the coefficients $A_j=0$ for $j\geq J+1$, where $J=p-q+n(q-1)$; see the appendix. As a consequence, the exponential series in (\ref{e38}) is then finite. In this case it is possible to evaluate the series for ${}_1\Psi_1(-x)$ in closed form and show that (\ref{e38}) becomes an {\em exact result\/}. To illustrate, we consider the particular case
$\alpha=\f{1}{3}$, $a=\f{1}{3}$ and $b=\f{7}{3}$, so that $\vartheta=-2$ and $J=2$. We obtain
\begin{eqnarray}
{}_1\Psi_1(-x)&=&\sum_{r=0}^\infty \frac{\g(\f{1}{3}r+\f{1}{3})}{\g(\f{1}{3}r+\f{7}{3})}\,\frac{(-x)^r}{r!}=9\sum_{r=0}^\infty (r+2)(r+3)\,\frac{(-x)^r}{(r+4)!}\nonumber\\
&=&\frac{9}{x^4}\left(\sum_{r=0}^\infty (r-1)(r-2)\,\frac{(-x)^r}{r!}+\frac{x^3}{3}-2\right)\nonumber\\
&=&3\left(\frac{1}{x}-\frac{6}{x^4}\right)+9x^{-2}e^{-x}\left(1+\frac{2}{x}+\frac{2}{x^2}\right)\label{e39}
\end{eqnarray}
upon evaluation of the infinite series as derivatives of $e^{-x}$. It is routine to verify that the right-hand side of (\ref{e39}) approaches the limit $\f{9}{4}$ as $x\ra 0$.
Insertion of the above parameter values into (\ref{e38}), together with $A_0=9$ and the normalised coefficients $c_j=A_j/A_0$ computed from the algorithm in the appendix $c_0=1$, $c_1=-2$, $c_2=2$ and $c_j=0$ ($j\geq 3$), yields the evaluation (\ref{e39}).

Finally, when $\vartheta=n$ it is seen from (\ref{e32}), or (\ref{e25}), that the algebraic expansion $H(x)\equiv 0$. In this case it is easily established that ${}_1\Psi_1(z)$ can be expressed as a polynomial in $z$ of degree $n$ multiplied by $e^z$. For, with $\Theta\equiv zd/dz$ and the polynomial $\wp_r(z)$ of degree $r\geq 1$ defined by
\[\wp_r(z):=e^{-z}\prod_{j=0}^{r-1}(\alpha\Theta+b+j)e^z,\]
we find from (\ref{e11}) 
\[{}_1\Psi_1(z)=\sum_{r=0}^\infty (\alpha r+b)_n\,\frac{z^r}{r!}
=\prod_{j=0}^{n-1}(\alpha\Theta+b+j)e^z=\wp_n(z)e^z.\]
Then, since $\wp_{r+1}(z) e^z=(\alpha\Theta+b+r)\wp_r(z) e^z$ and
for any differentiable function $\phi$ and constants $\alpha$ and $\gamma$
\[(\alpha\Theta+\gamma)\phi e^z=\{(\alpha z+\gamma)\phi+\alpha z\phi'\}\,e^z,\]
the polynomials $\wp_r(z)$ are defined recursively by
\[\wp_0(z)=1,\quad \wp_{r+1}(z)=(\alpha z+b+r)\wp_{r}(z)+\alpha z\wp_{r}'(z)\quad (0\leq r\leq n-1).\]

The algorithm described in the appendix shows that when $\vartheta=n$, the coefficients $A_j=0$ for $j\geq n+1$. It follows, from (\ref{e22c}), that the exponential series $E(z)$ is finite and by the uniqueness of polynomial representations we must have the exact result
\bee\label{e310}
{}_1\Psi_1\left(\begin{array}{c}a, \alpha \\  a-n, \alpha\end{array};z\right)=\wp_n(z)e^z=z^ne^z\sum_{j=0}^n A_jz^{-j}
\qquad(\vartheta=n)\ee
for $n=1, 2, \ldots\ $.

\vspace{0.6cm}

\begin{center}
{\bf 5. Numerical examples and concluding remarks}
\end{center}
\setcounter{section}{5}
\setcounter{equation}{0}
\renewcommand{\theequation}{\arabic{section}.\arabic{equation}}
We present some numerical examples to demonstrate the accuracy of the expansions in Theorems 3 and 4. We show in Table 1 the absolute relative error in the computation of
the exponentially small expansions in (\ref{e38a}) and (\ref{e38}) compared to their respective left-hand sides when $x=25$ for different parameter values and truncation index $j$. The value of ${}_1\Psi_1(-x)$ was obtained by high-precision computation of the series in (\ref{e11}). The optimal truncation index $m_o$ of the algebraic expansion was determined by inspection and it was verified that with this value of $x$ the range $j\leq 4$ corresponded to sub-optimal truncation of the exponentially small series.
\begin{table}[h]
\caption{\footnotesize{Values of the absolute relative error in the computation of the asymptotic exapnsion on the right-hand sides of (\ref{e38a}) and (\ref{e38}) when $x=25$ for different parameters and truncation index $j$. In the first two columns the algebraic expansion has been optimally truncated.}}
\begin{center}
\begin{tabular}{|c|c|c|c|}
\hline
&&&\\[-0.25cm]
\mcol{1}{|c|}{} & \mcol{1}{c|}{$a=\fs$,\,$b=\f{1}{4}$} & \mcol{1}{c|}{$a=\f{3}{4}$,\,$b=\f{1}{4}$}  & \mcol{1}{c|}{$a=\f{1}{2}$,\,$b=\f{5}{2}$}\\
&&&\\[-0.3cm]
\mcol{1}{|c|}{$j$} & \mcol{1}{c|}{$\alpha=\f{1}{3},\,\vartheta=+\f{1}{4}$} & \mcol{1}{c|}{$\alpha=\f{3}{2},\,\vartheta=+\fs$}  & \mcol{1}{c|}{$\alpha=\f{3}{4},\,\vartheta=-2$}\\
[.1cm]\hline
&&&\\[-0.3cm]
0  & $7.396\times 10^{-3}$ & $2.320\times 10^{-2}$  & $1.238\times 10^{-2}$ \\
1  &  $9.237\times 10^{-5}$ & $5.646\times 10^{-4}$& $1.116\times 10^{-3}$\\
2  &  $7.191\times 10^{-7}$& $3.738\times 10^{-5}$& $1.439\times 10^{-4}$\\
3  &  $6.365\times 10^{-7}$& $3.915\times 10^{-6}$& $2.413\times 10^{-5}$\\
4  &  $8.993\times 10^{-8}$& $5.533\times 10^{-7}$ & $4.987\times 10^{-6}$\\
[.2cm]\hline
\end{tabular}
\end{center}
\end{table}

The first column in Table 1 shows the relative error for $\vartheta=\f{1}{4}$, so that both series involving the coefficients $A_j$ and $B_j$ contribute, whereas the second column shows a case with $\vartheta=\fs$, so that only the series involving the coefficients $B_j$ contributes. The third column shows a case with $\vartheta=-2$ where the algebraic expansion consists of $n=2$ terms. The values of the coefficients $A_j$ used in these computations are presented in their normalised form in Table 2.

Two interrelated remarks can be made from Theorem 1. First, the Stokes multiplier (given by the quantity in curly braces in (\ref{e38a}) divided by $A_0$) is equal to $\cos \pi\vartheta$ to leading order. And secondly,
it is seen that the leading behaviour of the exponentially small expansion when $\kappa=1$ on the Stokes line $\arg\,z=\pi$  is O($x^\vartheta e^{-x})$ for $\vartheta$ not equal to half-integer values, but becomes O($x^{\vartheta-\fr}e^{-x})$ when $\vartheta=\pm\fs, \pm\f{3}{2}, \ldots\,$. 

Finally, when $\alpha=\beta=1$, the Wright function ${}_1\Psi_1(z)$ reduces to a multiple of the confluent hypergeometric function ${}_1F_1$ given by
\[{}_1\Psi_1\left(\begin{array}{c}a, 1\\b, 1\end{array};z\right)=\frac{\g(a)}{\g(b)}\,{}_1F_1(a;b;z).\]
Then, the expansion of ${}_1F_1(a;b;-x)$ for $x\ra+\infty$ can be obtained from Theorems 3 and 4, in which the coefficients $B_j$ are obtained by setting  $\mu=1$ in the coefficients $g_{2k}(\mu;j)$. It should be pointed out that the expansion in this case also follows from the exponentially improved expansion for the second Kummer function $U(a,b,z)$ combined with the connection formula (13.2.12) in \cite[p.~323, 329]{DLMF} and the expansion (\ref{e301}); see \cite{PCHF} for details.

\vspace{0.6cm}

\begin{center}
{\bf Appendix: \ An algorithm for the coefficients $c_j=A_j/A_0$}
\end{center}
\setcounter{section}{1}
\setcounter{equation}{0}
\renewcommand{\theequation}{\Alph{section}.\arabic{equation}}
We describe an algorithm for the computation of the normalised coefficients $c_j:=A_j/A_0$ appearing in the exponential expansion $E(z)$ in (\ref{e22c}).
The inverse factorial expansion (\ref{e22b}) can be written as
\bee\label{e41}
\frac{\g(s) \g(1-b+\beta s)}{\g(\kappa s+\vartheta) \g(1-a+\alpha s)}=\kappa A_0(h\kappa^\kappa)^{-s}\left\{\sum_{j=0}^{M-1}\frac{c_j}{(1-\kappa s-\vartheta)_j}+\frac{O(1)}{(1-\kappa s-\vartheta)_M}\right\}
\ee
for $|s|\ra\infty$ uniformly in $|\arg\,s|\leq\pi-\epsilon$. Introduction of the scaled gamma function $\g^*(z)=\g(z) (2\pi)^{-\fr}e^z z^{\fr-z}$ leads to the representation
\[\g(\alpha s+\gamma)=\g^*(\alpha s+\gamma) (2\pi)^\fr e^{-\alpha s} (\alpha s)^{\alpha s+\gamma-\fr} \,{\bf e}(\alpha s; \gamma),\]
where
\[{\bf e}(\alpha s; \gamma):= \exp\,\left[(\alpha s+\gamma-\fs) \log\,\left(1+\frac{\gamma}{\alpha s}\right)-a\right].\]

Then, with the definitions of the parameters in (\ref{e12}) and (\ref{e23}),
the above ratio of gamma functions becomes
\bee\label{e42}
\frac{\g(s) \g(1-b+\beta s)}{\g(\kappa s+\vartheta) \g(1-a+\alpha s)}=\kappa A_0(h\kappa^\kappa)^{-s}\,R(s)\,\Upsilon(s),
\ee
where
\[\Upsilon(s):=\frac{\g^*(s) \g^*(1-b+\beta s)}{\g^*(\kappa s+\vartheta) \g^*(1-a+\alpha s)},\qquad R(s):=\frac{{\bf e}(s; 0)\, {\bf e}(\beta s; 1-b)}{{\bf e}(\kappa s; \vartheta)\, {\bf e}(\alpha s; 1-a)}.\]
Substitution of (\ref{e42}) in (\ref{e41}) then yields
\bee\label{e43}
R(s)\,\Upsilon(s)=\sum_{j=0}^{M-1}\frac{c_j}{(1-\kappa s-\vartheta)_j}+\frac{O(1)}{(1-\kappa s-\vartheta)_M}
\ee
as $|s|\ra\infty$ in $|\arg\,s|\leq\pi-\epsilon$.

We now let $\chi:=(\kappa s)^{-1}$ and expand $R(s)$ and $\Upsilon(s)$ for $\chi\ra 0$ making use of the well-known expansion \cite[p.~71]{PK}
\[\g^*(z)\sim\sum_{k=0}^\infty(-)^k\gamma_kz^{-k}\qquad(|z|\ra\infty;\ |\arg\,z|\leq\pi-\epsilon),\]
where $\gamma_k$ are the Stirling coefficients, with 
\[\gamma_0=1,\quad \gamma_1=-\f{1}{12},\quad \gamma_2=\f{1}{288},\quad  \gamma_3=\f{139}{51840},
\quad \gamma_4=-\f{571}{2488320}, \ldots\ .\]
Expanding the right-hand side of (\ref{e43}) in powers of $\chi$, we can then match coefficients recursively with the aid of {\it Mathematica} to determine the $c_j$. This procedure is found to work well in specific cases when the various parameters have numerical values, where up to a maximum of 50 coefficients have been so calculated.
In Table 1 we present the normalised coefficients $c_j$ in the specific cases considered in Section 5.
\begin{table}[h]
\caption{\footnotesize{Values of the normalised coefficients $c_j=A_j/A_0$.}}
\begin{center}
\begin{tabular}{|c|l|l|l|}
\hline
&&&\\[-0.25cm]
\mcol{1}{|c|}{} & \mcol{1}{c|}{$a=\fs$,\,$b=\f{1}{4}$} & \mcol{1}{c|}{$a=\f{3}{4}$,\,$b=\f{1}{4}$}  & \mcol{1}{c|}{$a=\f{1}{2}$,\,$b=\f{5}{2}$}\\
&&&\\[-0.3cm]
\mcol{1}{|c|}{$j$} & \mcol{1}{c|}{$\alpha=\f{1}{3},\,\vartheta=+\f{1}{4}$} & \mcol{1}{c|}{$\alpha=\f{3}{2},\,\vartheta=+\fs$}  & \mcol{1}{c|}{$\alpha=\f{3}{4},\,\vartheta=-2$}\\
[.1cm]\hline
&&&\\[-0.3cm]
0  & \ \ \ 1  &  \ \  1 & \ \  1   \\
1  &  $-\f{3}{16}$ & $-\f{1}{8}$& $+\f{1}{3}$\\
[0.15cm]
2  &  $-\f{23}{512}$& $-\f{55}{1152}$& $+\f{7}{9}$\\
[0.15cm]
3  &  $+\f{343}{8192}$& $-\f{185}{3072}$& $+\f{70}{27}$\\
[0.15cm]
4  &  $+\f{133595}{524288}$& $-\f{351685}{2654208}$& $+\f{910}{81}$\\
[0.15cm]
5  &  $+\f{8169315}{8388608}$& $-\f{988855}{2359296}$& $+\f{14560}{243}$\\
[.2cm]\hline
\end{tabular}
\end{center}
\end{table}

Finally, when $\alpha=\beta$ ($\kappa=h=1$) and $\vartheta=n$ ($n=1, 2, \ldots$)  the quotient of gamma functions in (\ref{e41}) reduces to
\[A_0\prod_{j=1}^n \frac{s-(b+j-1)/\alpha}{s+j-1}.\]
It is seen by comparison with the right-hand side of (\ref{e41}) that the series involving the coefficients $c_j$ must terminate when $j=n$; that is, $c_j=0$ for $j\geq n+1$.
In a similar manner, when $\vartheta=-n$ and the parameters $\alpha$, $a$ satisfy (\ref{e4c}), with $p$, $q$ denoting relatively prime positive integers, the quotient  in (\ref{e41}) reduces to
\[A_0\prod_{j=1}^n\frac{s-j}{s-p-(j-1)q}.\]
In this case, the series on the right-hand side of (\ref{e41}) terminates when $j=J$, where $J=p-q+n(q-1)$; that is, $c_j=0$ for $j\geq J+1$.
The coefficients $c_j$ can be determined by recursive application of the `cover-up' rule in partial fractions; in particular, we have 
\[c_n=\alpha^{-n} (b)_n\quad (\vartheta=n),\qquad c_J=(-)^J \alpha^{n-1}(n)_J\quad(\vartheta=-n).\].

\end{document}